\documentclass[letterpaper,10pt]{amsart}
\usepackage{amssymb}

\newcommand{\mb}[1]{\mbox{\rm {#1}}}
\newcommand{\mr}[1]{\mathfrak{#1}}

\newcommand{\wId}[1]{\equiv_w 1 (\mb{mod } {#1})}
\newcommand{\wsId}[1]{\equiv_w^* 1 (\mb{mod } {#1})}
\newcommand{\wre}[1]{\preceq_{#1}}

\newcommand{\Tsub}{\unlhd_T}

\newcommand{\Z}{\mathbb{Z}}
\newcommand{\N}{\mathbb{N}}

\newtheorem{lm}{Lemma}[section]
\newtheorem{thq}[lm]{Theorem}
\theoremstyle{remark}
\newtheorem{rem}[lm]{Remark}
\newtheorem{ex}[lm]{Example}
\newtheorem{df}[lm]{Definition}
\newtheorem{co}[lm]{Corollary}
\newtheorem{que}[lm]{Open Question}

\DeclareMathOperator{\id}{id}

\DeclareMathOperator{\Cen}{Cen}

\newcommand{\la}{\langle}
\newcommand{\ra}{\rangle}

\title{Weak Identities in finitely generated groups}
\author{Martin Kassabov}
\address{
Department of Mathematics \\
University of Alberta\\
632 Central Academic Building\\
Edmonton, Alberta, T6G 2G1\\
Canada}

\email{kassabov@aya.yale.edu}
\urladdr{http://www.math.ualberta.ca/\~{}mkassabov}

\begin{document}


\begin{abstract}
In this article we introduce the notion of weak identities in a
group and study their properties. We show that weak identities
have some similar properties to ordinary ones. We use this notion to
prove that any finitely generated solvable discriminating group is
abelian, which answers a question raised in~\cite{Disc2}.
\end{abstract}

\maketitle

\begin{section}{Introduction}

In this article we introduce the notion of weak identities in a
group and study their properties. An element $f(g_1,\dots, g_k)$
in the free group on $k$ generators is called a weak identity in group $G$
if there exists $N$, such that for any $N$ $k$-tuples
$(h_{i1}, \dots ,h_{ik})_{i=1,N}$ of elements in the group $G$ such that
any two elements in different $k$-tuples commute, then $f$
evaluated at one $k$-tuple gives identity\footnote{The main
differences between the weak identities and the ordinary ones
come form the fact that it is only known that $f$ vanishes on
some $k$-tuple but it is not known on which one.}.

First, we show that the set of weak identities in a given group $G$ form a
verbal subgroup. This result shows that weak identities are similar
to the ordinary ones, but there are some substantial differences -- the main ones are:
\begin{list}{$\bullet$}{}
\item
the free group has many nontrivial weak identities;

\item
the class of groups satisfying a given weak identity is not closed under
taking homomorphic images.
\end{list}

We show that a very large class of groups satisfy the weak
identity $[g_1,g_2]$ -- this class includes all linear groups.
Next, we study weak identities modulo a verbal subgroup, and use
them to construct a relation on all verbal subgroups of a free
group. We also introduce weak* identities, since the above
relation is not transitive. Finally, we use the notion of weak
identity to study discriminating groups and answer two questions
raised in~\cite{Disc2}.

The paper is organized as follows:
The notion of weak identities is defined in section~\ref{sec:weak}.
In section~\ref{sec:bcs}, we discuss the notion of a group
having a bounded centralizer sequence and show the connection with
weak identities.
In section~\ref{sec:weak*}, we showed that the notion of weak identities is
not transitive and define weak* identities in order to address
this problem.
In section~\ref{sec:meta}, we investigate weak* identities in finitely
generated meta abelian groups.
In section \ref{sec:disc}, we apply the results from the previous sections
to discriminating groups and give answers
to the Questions 2D and 3D from \cite{Disc2} - we prove that:
every linear discriminating group is abelian;
and that every  finitely generated solvable discriminating  group is
free abelian.
Finally, in section~\ref{sec:open} we pose some
open questions concerning the notion of weak identities.

Acknowledgments: The author thanks R. Muchnik for
the useful discussions and A. Myashnikov for
introducing the author to the topic.

\end{section}


\begin{section}{Weak Identities}
\label{sec:weak}

In this section we define the notion of weak identities in a group.
They have similar properties to ordinary identities -- see theorems~\ref{gen}
and~\ref{composition}; however there are also few substantial differences --
see remarks~\ref{widentities} and~\ref{wvariety}.

\begin{df}
Let $\mr{F}$ be the free group on countably many generators $g_i$,
for $i\in \N$.
A subgroup $\mr{H}$ in $\mr{F}$ is called a $T$-subgroup or verbal subgroup
(denoted $\mr{H} \Tsub \mr{F}$)
if it is
preserved by all endomorphisms of the group $\mr{F}$.
For a set $S \subset \mr{F}$, we will denote by $\la S \ra_T$
the minimal $T$-subgroup which contains $S$.

\end{df}

\begin{df}
Let $G$ be an abstract group. We say that the set $S \subset F$ is
a set of weak identities in the group $G$, if there exists an integer
$N$ such that for any
elements $s_k \in S$, for $k=1,\dots, N$, and any homomorphism
$$
\rho: \mr{F}^{\times N} = \mr{F} \times \mr{F} \times \dots \times \mr{F}
\to G,
$$
there exists an index $k$ between $1$ and $N$ such that
$$
\rho(i_k(s_k)) = 1,
$$
where $i_k$ denotes the inclusion of $\mr{F}$ in $\mr{F}^{\times N}$
at the $k$-th component.
The number $N$ is called the height of the set $S$ of weak identities.

An element $f\in \mr{F}$ is called a weak identity in $G$,
if the set $\{ f \}$ is a set of weak identities.
\end{df}

\begin{co}
Let $S_i$ be a finite collection of sets, such that $S_i$ is a set of weak
identities in $G$ of height $n_i$, for each $i$.
Then the union $\cup S_i$ is also a set of weak identities of height
at most $\sum n_i$.
\end{co}

\begin{ex}
\label{free} The element $[g_1,g_2]= g_1g_2 g_1^{-1} g_2^{-1}$ is
a weak identity of height $2$ in the free group $\mr{F}$, but it
is not an ordinary identity. In order to prove this we need to
show that for any homomorphism  $\rho : \mr{F}^{\times 2} \to
\mr{F}$, we have that $\rho(i_1[g_1,g_2]) =1$ or
$\rho(i_2[g_1,g_2]) =1$. Let us assume that $\rho(i_1([g_1,g_2]))
\not= 1$. Then the centralizer of this element in the free group
$\mr{F}$ is an infinite abelian group $H$. The elements
$\rho(i_2(g_1))$ and $\rho(i_2(g_2))$ lie in $H$ because they
commute with $\rho(i_1([g_1,g_2]))$. This implies that their
commutator is $1$. The above argument shows that at least one of
the elements $\rho(i_1([g_1,g_2]))$ or $\rho(i_2([g_1,g_2]))$ is
$1$, i.e., $[g_1,g_2]$ is a weak identity in the free group
$\mr{F}$ of height $2$. This example will be generalized in
section~\ref{linear}.
\end{ex}

The next theorem states that in order to show that a $T$-subgroup
consists of weak identities, it is enough to verify that the set of
its generators are form a set of weak identities.
\begin{thq}
\label{gen}
If  $S$ is a set of weak identities in the group $G$, then the $T$-subgroup
$\mr{S}= \la S \ra_T \Tsub \mr{F}$ generated by the set $S$ is also
a set of weak identities.
\end{thq}
\begin{proof}
Any element $g$ in the $T$-subgroup $\mr{S}$, generated by $S$,
can be written as the product
$$
g = \prod_{k=i}^n \phi_i(s_i),
$$
where $s_i\in S$ and $\phi_i: \mr{F} \to \mr{F}$ are endomorphisms.

Suppose that we have $N$ elements $g_j \in \mr{S}$ and a homomorphism
$\rho: \mr{F}^{\times N} \to G$. We can write any of the elements
$g_j$ in the form
$$
g_j= \prod_{i=1}^{n_j} \phi_{i,j}(s_{i,j}).
$$
Let $i(j)$ be the index between $1$ and $n_j$ such that the length of
$\rho(i_j(\phi_{i,j}(s_{i,j})))$ is maximal.\footnote{
In order to have a notion of length of an element in $G$, we need to fix
a generating set of $G$.}

Now consider the homomorphism $\tilde \rho : \mr{F}^{\times N} \to G$
defined by
$$
\tilde \rho = \rho \circ ( \phi_{i(1),1} \times \dots \times \phi_{i(N),N}).
$$
By definition we have that
$\tilde \rho(i_j (s_{i(j),j})) = \rho(i_j(\phi_{i(j),j}(s_{i(j),j})))$.

Using the homomorphism $\tilde \rho$ and fact that $S$ is a set of
weak identities in $G$ (by construction we have $s_{i(j),j} \in
S$, for all $j$), we know that there exists $j$ such that $\tilde
\rho(i_j (s_{i(j),j}))=1$.

By the definition of $i(j)$, we have that the length of
$\rho(i_j(\phi_{i(j),j} (s_{i(j),j}))$ is bigger than or equal to
the length of $\rho(i_j (\phi_{i,j}(s_{i,j})))$ for any $i$.
However, the first element is identity and has length zero,
therefore the lengths of all elements $\rho(i_j
(\phi_{i,j}(s_{i,j})))$ are $0$, i.e., all of hem are equal to the
identity in $G$. This shows that
$$
\rho(i_j(g_j)) = \prod \rho(i_j (\phi_{i,j}(s_{i,j}))) = 1,
$$
which shows that the subgroup $\mr{S}$ is a set of weak identities in $G$.
\end{proof}

\begin{rem}
\label{widentities}
Let us fix a group $G$. Denote by $\mr{wId}(G)$ the set of all
elements $f$ in $\mr{F}$,
such that the set $\{f\}$ is a set of weak identities in $G$.
By the previous theorem, this is a $T$-subgroup in $\mr{F}$ which is
called the group of weak identities in $G$. Note that this theorem does
not imply that $\mr{wId}(G)$ is a set of weak identities in $G$. However,
any finitely generated\footnote{as T-subgroup} $T$-subgroup
$\mr{H}$ of $\mr{wId}(G)$
is a set of weak identities.
\end{rem}

\begin{ex}
\label{finite}
If the group $G$ is finite, then $\mr{wId}(G) = \mr{F}' . \mr{F}^n$, where
$n$ is the minimal number such that $g^n = 1$ for any $g \in G$.

First, let us show that $[g_1,g_2]$ is a weak identity in $G$.
Suppose that it is not. Then for any $N$, there exist elements
$g_i$ and $h_i$, for $1\leq i \leq N$, in the group $G$, such that
$$
[g_i,g_j]=1 \quad [h_i,h_j]=1 \quad [g_i,h_j]=1, \mb{ iff }i\not=j
\quad [g_i,h_i] \not= 1.
$$
Let us define the subgroups $P_i$ of $G$ using these elements
$g_i$ by
$$
P_i = \{ g\in G | [g,g_j]=1 \mb{ for }j<i \},
$$
It is easy to check that $h_i \in P_{i-1} \setminus P_i$, i.e.
$P_i$ form a strictly descending sequence of subgroups in $G$,
which is impossible if $N > \log_2 |G|$. This contradiction shows
that $[g_1,g_2]$ is a weak identity in any finite group $G$ (of
height at most $\log_2 |G|$).

Since the element $g^n$ is an identity in $G$, it is also a weak
identity (of height $1$). Therefore, the group  $\mr{wId}(G)$
contains the subgroup $\mr{F}' . \mr{F}^n$. If we assume that the
inclusion is strict, then the group $\mr{wId}(G)$ would contain
the element $g^k$ for some $0<k<n$. The last element is not an
identity in $G$. It can be shown that it is not also a weak
identity. This proves that $\mr{wId}(G) = \mr{F}' . \mr{F}^n$.
\end{ex}

\begin{rem}
\label{wvariety} Let $\mr{H}$ be a $T$-subgroup. Denote by
$\mr{wVar}(\mr{H})$ the class of all groups $G$ such that any
element in $\mr{H}$ is a weak identity in $G$\footnote{This is the
same as saying that any finite subset is a set of week identities
in the group $G$.}. It can be shown that the class $\mr{wVar}(\mr{H})$ is
closed under taking subgroups
and taking finite Cartesian products, 
but not infinite products. In general is not closed under taking
homomorphism images, although it is closed under taking some kinds
of restricted homomorphic images (see open
problem~\ref{describe}). This is one important difference between
weak identities and the ordinary ones, since it implies that there
are no universal objects in the class $\mr{wVar}(\mr{H})$.
\end{rem}

\end{section}


\begin{section}{Weak identities in linear groups}
\label{linear}
\label{sec:bcs}

In this section we show that a large class of groups
lie in $\mr{wVar}(\mr{F}')$. This generalizes
examples~\ref{free} and~\ref{finite}.

\begin{df}
A group $G$ is said to have bounded centralizer sequences if there
exists an integer $N$, such that any strictly increasing sequence
of stabilizers of sets of mutually commuting elements has length
less than $N$. That is, for any sequence of subsets
$\{P_i\}_{i=1}^N$ of $G$ such that
$$
P_1 \subset P_2 \subset \dots \subset P_N,
$$
and $[P_i,P_i]=1$ for any $i$, we have that the sequence of their
centralizers
$$
\Cen(P_1) \supset \Cen(P_2) \supset \dots \supset \Cen(P_N)
$$
is not strictly decreasing.
\end{df}

\begin{ex}
Any finite group has bounded centralizer sequences and the same is true for
any free group. 
\end{ex}

\begin{lm}
\label{bcs}
If the group $G$ has bounded centralizer sequences, then $[g_1,g_2]$
is a weak identity in $G$.
\end{lm}
\begin{proof}
Suppose that $[g_1,g_2]$ is not a weak identity in $G$.
Then for any $N$ there exist elements $g_i$ and $h_i$, for $1\leq i \leq N$
in the group $G$, such that
$$
[g_i,g_j]=1, \quad [h_i,h_j]=1, \quad [g_i,h_j]=1 \mb{ iff }i\not=j,
\quad [g_i,h_i] \not= 1.
$$
The centralizers of the sets
$P_i = \{ g_1,\dots,g_i\}$ of mutually commuting elements, form a
strictly descending sequence, because
$h_i \in \Cen(P_{i-1}) \setminus \Cen(P_i)$. This contradicts
the assumption that the group $G$ has bounded centralizer sequences.
Therefore, $[g_1,g_2]$ is a weak identity in the group $G$
of height equal to the maximal length of an increasing
sequence of centralizers in $G$.
\end{proof}

\begin{thq}
\label{lin}
If $G$ is a linear group, then $\mr{F}'$ consists of weak identities in $G$.
\end{thq}
\begin{proof}
By lemma~\ref{bcs} it is enough to show that any linear group has
a bounded sequence of centralizers.
The centralizer of a subset $P$ of the linear group $GL_n$
is the same as the centralizer of the linear
span of $P$ in the matrix algebra $M_n$.
But in the matrix algebra $M_n$ any strictly increasing sequence of sub-spaces
has bounded length. Therefore, any strictly decreasing sequence of
centralizers has a bounded length.
\end{proof}
\end{section}


\begin{section}{Weak* Identities}
\label{sec:weak*}

It is also possible to define weak identities modulo some
$T$-subgroup $\mr{H}$. In order to do this, we need to define
the normal subgroup $\mr{H}(G)$ of $G$.

\begin{df}
Let $\mr{H} \Tsub \mr{F}$ be a $T$-subgroup, and let $G$ be a
group. Denote by
$$
\mr{H}(G) = \{ \pi(h) \mid h\in \mr{H},\,\, \pi : \mr{F} \to G \}
$$
the subgroup of $G$ consisting of all elements which are the
images of elements in $\mr{H}$ under homomorphisms from $\mr{F}$
to $G$.
\end{df}

Let us fix a finitely generated group $G$. For the rest of this
section all the identities we consider are in the group $G$,
unless stated otherwise.

\begin{df}
A set $S$ is said to be a set of weak identities (in the group
$G$) modulo the $T$-subgroup $\mr{H}$, denoted $S \wId{\mr{H}}$,
if $S$ is a set of weak identities in the group $G/\mr{H}(G)$.
This gives rise to a relation on verbal subgroups of the free
group: we say that $\mr{H}_1 \wre{wG} \mr{H}_2$ iff $\mr{H}_1
\wId{\mr{H}_2}$ in the group $G$. In the case of ordinary
identities the above relation comes from the inclusion of the
corresponding subgroups, i.e., $\mr{H}_1 \wre{G} \mr{H}_2$ iff
$\mr{H}_1(G) \subseteq \mr{H}_2(G)$.
\end{df}

\begin{rem}
Note that if one takes the class of groups $G$ such that $S$ is a
set of weak identities modulo $\mr{H}$, this class is not closed
under taking subgroups. The same is true if one considers ordinary
identities.
\end{rem}

In order to state Theorem~\ref{composition}, we need to describe
one
construction of verbal subgroups.

\begin{df}
Let $\mr{F}_n$ be the free group generated by $g_k$, for
$k=1,\dots,n$, and let $f\in \mr{F}_n$ be an element in it. For a
$T$-subgroup $\mr{H}$ and an index $1\leq i \leq n$, denote by
$\la f_{|g_i\to \mr{H}} \ra_T$ the $T$-subgroup generated by the
set
\begin{equation}
\label{sub} \{ \rho(f) \mid \rho: \mr{F}_n \to \mr{F},
\rho(g_j)=g_j, \mb{ for }j\not=i, \,\,\rho(g_i) \in \mr{H} \},
\end{equation}
i.e., all the elements which can be obtained from the element $f$
by substituting a word from $\mr{H}$ in the place of $g_i$.

Similarly, if $\mr{H}_i$ are $T$-subgroups, for $i=1,\dots,n$,
then by $\la f_{|g_i\to \mr{H}_i} \ra_T$ we will denote the
$T$-subgroup in $\mr{F}$ generated by
\begin{equation}
\label{submany} \{ \rho(f) \mid \rho: \mr{F}_n \to \mr{F},\,\,
\rho(g_i) \in \mr{H_i} \mb{ for all } i\}.
\end{equation}
\end{df}


\begin{thq}
\label{composition} a) Let $f\in \mr{F}_n$ and  let $\mr{S}$ and
$\mr{H}$ be $T$-subgroups in $\mr{F}$. If $\mr{S} \wId{\mr{H}}$,
then $\la f_{|g_i\to \mr{S}} \ra_T \wId{\la f_{|g_i\to \mr{H}}
\ra_T}$.

b) Let $f\in \mr{F}_n$ and  let $\mr{S}_i$ and $\mr{H_i}$ be
$T$-subgroups in $\mr{F}$. If $\mr{S}_i \wId{\mr{H}_i}$ for every
$i$, then $\la f_{|g_i\to \mr{S}_i} \ra_T \wId{\la f_{|g_i\to
\mr{H}_i} \ra_T}$.
\end{thq}
\begin{proof}
a) Let $N$ be the height of the set $\mr{S}$ of weak identities
modulo $\mr{H}$. First, we will show that the set~(\ref{sub}) of
generators of the $T$-subgroup $\la f_{|g_i\to \mr{S}} \ra_T$ form
a set of weak identities modulo $\la f_{|g_i\to \mr{H}} \ra_T$ of
height $N$. Suppose that
$$
a_k=\pi_k(f), \mbox{ where }\pi_k: \mr{F}_n \to
\mr{F},\,\,\pi_k(g_j)=g_j, \mb{ for }j\not=i, \,\,\pi_k(g_i) \in
\mr{S},
$$
are elements of type~(\ref{sub}) and that $\rho:\mr{F}^{\times N}
\to G$ is a homomorphism. Since $\mr{S} \wId{\mr{H}}$ and
$\pi_j(g_i) \in \mr{S}$, there exists an index $j$ such that
$\rho(i_j(\pi_j(g_i))) \in \mr{H}(G)$. Therefore, there exists an
element $h\in \mr{H}$ and a homomorphism $\bar \pi: \mr{F} \to G$
such that
$$
\rho(i_j(\pi_j(g_i))) = \bar \pi (h).
$$
Without loss of generality, we may assume that the element $h$
does not depend on the letters $g_k$ for $k\leq n$. Let us define
$\tilde \pi: \mr{F}_n \to \mr{F}$ to be a homomorphism which sends
$g_k$ to $g_k$ for $k\not=i$ and $\tilde \pi(g_i) = h$. Also
define $\tilde{\rho}: \mr{F} \to G$ by
$$
\tilde{\rho}(g_i) = \rho(i_j(g_i)), \mb{ for } i\leq n \quad
\mb{and} \quad \tilde{\rho}(g_i) = \bar\pi(g_i),\mb{ for } i>n.
$$
Then we have that
$$
\tilde \rho (\tilde \pi(f)) = \rho(i_j(a_j)).
$$
However, $\tilde \pi (f) \in \la f_{|g_i\to \mr{H}} \ra_T$ because
$\tilde \pi (g_i) \in \mr{H}$. This shows that
$$
\rho(i_j(a_j)) \in \la f_{|g_i\to \mr{H}} \ra_T (G).
$$
This proves that the generators of the $T$-subgroup $\la
f_{|g_i\to \mr{S}} \ra_T$ are weak identities of height $N$ in
$G/\la f_{|g_i\to \mr{H}} \ra_T(G)$. Finally, we can use
lemma~\ref{gen} to show that the whole $T$-subgroup consist of
weak identities.

b) Let $N_i$ be the heights of the sets $\mr{S_i}$ as weak
identities modulo $\mr{H_i}$. Suppose that
$$
a_k=\pi_k(f), \mbox{ where }\pi_k: \mr{F}_n \to \mr{F},
\,\,\pi_k(g_i) \in \mr{S_i},
$$
are elements of type~(\ref{submany}), where $N = \sum N_i$. Also
suppose that $\rho:\mr{F}^{\times N} \to G$ is a homomorphism.
Using that $\mr{S_i} \wId{\mr{H_i}}$ and $\pi_j(g_i) \in \mr{S_i}$
and that $N$ is big enough, we can show that there exists an index
$j$ such that $\rho(i_j(\pi_j(g_i))) \in \mr{H_i}(G)$, for every
$i$. Now we can repeat the proof of part a) to show that
$$
\rho(i_j(a_j)) \in \la f_{g_i\to \mr{H}_i} \ra_T (G).
$$
\end{proof}

There is one substantial difference between weak identities and ordinary ones:
if we have three verbal subgroups $\mr{H}_1$, $\mr{H}_2$ and $\mr{H}_3$, such
that $\mr{H}_1$ consists of weak identities modulo $\mr{H}_2$ and
$\mr{H}_2$ consists of weak identities modulo $\mr{H}_3$, then it does not
follow that $\mr{H}_1$ consists of weak identities modulo $\mr{H}_3$.

\begin{ex}
Let $G$ be nonabelian finite simple group. Then as shown in example~\ref{finite},
the element $[g_1,g_2]$ is a weak identity in $G$ (modulo the trivial verbal subgroup).
But the whole free group $\mr{F}$ consists of weak identities modulo
$\mr{F}'$ because the quotient $G/[G,G]$ is the trivial group. Notice that
$\mr{F}$ does not consist of weak identities in $G$ because the group $G$ is not
trivial. This example shows that in general the relation $\mr{H}_1 \wre{wG} \mr{H}_2$
is not transitive.
\end{ex}

In order to address this problem we need to define weak* identities.

\begin{df}
Let  $\mr{S}$ and $\mr{H}$ be $T$-subgroups. We call $\mr{S}$ weak* identities
modulo $\mr{H}$, denoted
$\mr{S} \wsId {\mr{H}}$, if there exist integer $n$ and  $T$-subgroups
$\mr{S}_i$, for
$i=0,\dots, n$ such that
$\mr{S} = \mr{S}_0$, $\mr{H}=\mr{S}_n$,
and $\mr{S}_{i-1}$ consists of weak identities modulo $\mr{S}_i$,
for all $i=1,\dots,n$.
The integer $n$ is called the length of the $T$-subgroup $\mr{S}$ as weak* identities
modulo $\mr{H}$.

We say that $\mr{S}$ consists of weak* identities in the group $G$
iff $\mr{S} \wsId{ \{1\}}$, i.e., if they are identities modulo the trivial
verbal subgroup.
\end{df}


\begin{thq}
\label{properties}
a)
Let $\mr{H_i} \Tsub \mr{F}$, for $i=1,\dots ,n$ be a descending
chain of $T$-subgroups, i.e.,
$$
\mr{H}_1 \supset \mr{H}_2 \supset \dots \supset \mr{H}_{n-1} \supset \mr{H}_n.
$$
Then
$\mr{H}_1 \wsId{\mr{H}_k}$ if and only
if $\mr{H}_i \wsId{\mr{H}_{i+1}}$ for every $i=1,\dots ,n-1$.

b)
Let  $f\in \mr{F}_n$ be a word. Then
$\la f_{|g_i\to \mr{S}} \ra_T$ consists of weak* identities modulo
$\la f_{|g_i\to \mr{H}} \ra_T$,
provided that $\mr{S} \wsId{\mr{H}}$.

c)
Let  $f\in \mr{F}_n$ be a word and let $\mr{S}_i$ and $\mr{H_i}$
be $T$-subgroups
in $\mr{F}$.
Then $\la f_{|g_i\to \mr{S}_i} \ra_T$ consists of weak* identities modulo
$\la f_{|g_i\to \mr{H}_i} \ra_T$,
whenever $\mr{S}_i \wsId{\mr{H}_i}$ for every $i$.

\end{thq}
\begin{proof}
a)
The ``if'' part follows from the definition of weak* identities;
the ``only if'' part  is trivial.

b), c) Induction on the length of $\mr{S}$ as a weak* identify modulo
$\mr{H}$ (the maximal of the length for part c) ) - the base case is trivial,
and the induction step follows theorem~\ref{composition}.
\end{proof}

\end{section}


\begin{section}{Weak* identities in solvable groups}
\label{sec:meta}

The main results in this section are theorems~\ref{meta} and~\ref{main},
which show that in
meta-abelian/solvable groups there are many weak* identities which are not
ordinary identities.

\begin{thq}
\label{meta}
If the group $G$ is finitely generated, then $\mr{F}'$ consists of
weak* identities in $G$
modulo $\mr{F}''$.
\end{thq}
\begin{rem}
The above theorem is not true for an infinitely generated group.
Let $G$ be any group.
Then any weak/weak* identity in the group $G^{\times \infty}$ is an
identity in $G$
and in $G^{\times \infty}$ (this is true because the group $G^{\times \infty}$
is discriminating and we can apply the results from the next section).
Therefore, there exists an infinitely generated
meta-abelian group $\widetilde G$ such that $[g_1,g_2]$ is not a weak* identity in
$\widetilde G$.
\end{rem}
\begin{proof}
In order to prove the theorem we need to prove several lemmas.
Let us fix the finitely generated group
$G$, generated by the $d$ elements $p_1,\dots,p_d$.
Without loss of generality, we may assume that $G$ is a meta-abelian group.

\begin{lm}
\label{nil3.5}
There exists an integer $N$, depending on the group $G$, such that the element
$[[g_1,g_2],g_3^N]$ is a weak identity in $G$ (modulo $\mr{F}''$).
\end{lm}
\begin{proof}
The group $G/G'$ acts by conjugation on $G'$, and $G'$ is a finitely generated
module over it.
Therefore, there exists an integer $N$, such that
for every $\bar g \in G/G'$ and $h \in G'$,
the length of the orbit of the element $h$ under the action of $\bar g$ is
either infinity or divides $N$.
Thus, for every $g \in G$ and $h \in G'$, $[h,g^N] \not =1$ implies that
$[h,g^k] \not =1$, for all $k\not=0$.

Suppose that we have a homomorphism $\rho$ from $\mr{F}^{\times d+1}$ to $G$,
such that
$$
\rho(i_k([[g_1,g_2],g_3^N])) \not= 1, \, \mb{ for all } k.
$$
Let $a_k=\rho(i_k(g_3))$ and
$b_k=\rho(i_K([g_1,g_2]))$. We have that $[b_k,a_j^N]\not = 1$ iff $k=j$.
From the choice of $N$ we have that $[b_k,a_k^n]\not = 1$ for all $n$.
This shows that the element $a_k^n$ does not lie in the subgroup generated by
$a_j$, for $j\not=k$ and $G'$,
because all elements in this subgroup commute with $b_k$.
Therefore,
the $a_i$-es are linearly independent modulo $G'$, which is impossible,
since the group $G/G'$ is an abelian group generated by $d$ elements.
Therefore, there is no such
homomorphism $\rho$, i.e., $[[g_1,g_2],g_3^N]$ is a weak identity.
\end{proof}

\begin{lm}
\label{nil3}
The element
$[[g_1,g_2],g_3]$ is a weak identity modulo the $T$-subgroup
$\mr{H}_3$ generated by
$\mr{F}''$ and $[[g_1,g_2],g_3^N]$.
\end{lm}
\begin{proof}
Let $\rho: \mr{F}^{\times n} \to G$ be a homomorphism, where $n = d\log_2 N$.
Suppose that
$\rho(i_k([[g_1,g_2],g_3])) \not \in \mr{H}_3(G)$ for all $k$. Let
$$
H_k = \{g\in G | [g,\rho(i_j([g_1,g_2]))] \in \mr{H}_3(G), \,
\mb{ for }j\leq k\}.
$$
The groups $H_k$ form a descending chain of subgroups.
The chain is strictly descending, since
$\rho(i_k(g_3)) \in H_{k-1} \backslash H_k$. All these subgroups
contain the subgroup
generated by $G'$ and the elements $g^N$, for any $g\in G$.
Therefore, we can project the subgroups $H_k$ to subgroups of the
group $G/G'.G^N$.
The last group is a finite group containing less than $N^d$ elements
and does not have a strictly
decreasing sequences of subgroups of length more than $d\log_2 N$ --
contradiction.
This
proves that the element $[[g_1,g_2],g_3]$ is a weak identity modulo $\mr{H}_3$.
\end{proof}

\begin{lm}
\label{nil2.5}
There exists an integer $N$, depending on the group $G$, such that the element
$[g_1,g_2^N]$ is a weak identity modulo the $T$-subgroup $\mr{H}_2$ generated
by $[[g_1,g_2],g_3]$.
\end{lm}
\begin{proof}
The proof is similar to the one of lemma~\ref{nil3.5}.
The group $G'/\mr{H}_2$ is a finitely generated abelian group.
Let $N$ be a number divisible by the
order of all torsion elements in this group. This implies that if $
h \in G'$ and
$h^N \not \in \mr{H}_2(G)$, then $h^k \not \in \mr{H}_2(G)$, for all $k\not=0$.
Now if we assume that $[g_1,g_2^N]$ is not a weak identity, then there
exists a map
$\rho: \mr{F}^{\times d+1} \to G$ such that
$\rho(i_k([g_1,g_2^N])) \not \in \mr{H}_2(G)$.
Using arguments similar to the ones in the proof of lemma~\ref{nil3.5},
we can show
that the elements $\rho(i_k(g_2))$
are linearly independent modulo $G'$, which is impossible -- contradiction.
\end{proof}

\begin{lm}
\label{nil2}
The element
$[g_1,g_2]$ is a weak identity modulo the $T$-subgroup $\mr{H}_1$ generated
by $[[g_1,g_2],g_3]$ and $[g_1,g_2^N]$
\end{lm}
\begin{proof}
Same as the proof of lemma~\ref{nil3}, but we use the element $g_1$,
instead of the element $[g_1,g_2]$ to construct the groups $H_k$.
\end{proof}

\begin{rem}
The statements of Lemmas~\ref{nil3} and~\ref{nil2} hold in any
finitely generated group $G$, although the heights of the
corresponding weak identities depend only on the number of
generators of the group $G$. The same is not true for
Lemmas~\ref{nil3.5} and~\ref{nil2.5}, where the number $N$ depends
on the group $G$.
\end{rem}

Now we can use these lemmas to prove Theorem~\ref{meta}: take the $T$-subgroups
$$
\mr{F}' \supset \mr{H}_1 \supset \mr{H}_2 \supset \mr{H}_3 \supset \mr{F}''.
$$
Lemmas~\ref{nil3.5}, \ref{nil3}, \ref{nil2.5} and \ref{nil2}
show that each of these $T$-subgroups consists
of weak identities modulo the next one. Therefore, by definition
$\mr{F}'$ consists of
weak* identities modulo $\mr{F}''$, which completes the proof of
theorem~\ref{meta}.
\end{proof}

\begin{co}
\label{solvable}
If for a fixed finitely generated group $G$,
the $T$-subgroup $\mr{F}'$ consists of
weak* identities modulo $\mr{F}''$.
Then $\mr{F}'$ consist of weak* identities modulo ${\mr{F}^n}'$ for all $n$.
Here, $\mr{F}'$ denotes the commutator subgroup,
$\mr{F}''$ is its commutator subgroup and ${\mr{F}'}^n$ is the $n$-th term
of the solvable series of the free group $\mr{F}$.
\end{co}
\begin{proof}
We can use theorem~\ref{properties} c) to show that  ${\mr{F}'}^k$ consists of
weak* identities in $G$ modulo the group ${\mr{F}'}^{k+1}$. Finally,
using part a)
we can show that
$$
\mr{F}' \wsId{{\mr{F}'}^n}.
$$
\end{proof}

\begin{thq}
\label{main}
Let $G$ be a finitely generated solvable group. Then $\mr{F}'$ consists
of weak* identities.
\end{thq}
\begin{proof}
Corollary~\ref{solvable} and theorem~\ref{meta} imply that
$\mr{F}' \wsId{\mr{F}'}^n$. The group $G$ is solvable. Therefore,
${\mr{F}'}^n(G)=1$, and all identities mod ${\mr{F}'}^n$ are identities in $G$.
\end{proof}
\end{section}


\begin{section}{Discriminating groups}
\label{sec:disc}

Discriminating groups were introduced in~\cite{Disc1} and~\cite{AG} by
G. Baumslag, A. Myasnikov and  V. Remeslenikov.
There the authors consider the universal theory of groups. A group $G$
is called square-like if the universal theories of $G$ and $G \times G$
coincide. There is simple sufficient condition of a group $G$ to be
square-like - a group satisfying this condition is called discriminating.
More detailed discussion
of discriminating groups can be found in~\cite{AG}, \cite{Disc1} and~\cite{Disc2}.

A group $G$ discriminates a group $H$ if for any finite set of
nontrivial elements $h_i \in H$ there exists a
homomorphism $\phi:H \rightarrow G$ which maps $h_i$-es to nontrivial elements in $G$.
A group $G$ is called discriminating, if $G$ discriminates $G \times G$.
The definition we give below is not exactly the same as the
one in~\cite{AG}, but it is equivalent.

\begin{df}
A group $G$ is called discriminating if for any integer $N$ and for
any elements
$h_1,\dots,h_N \in G \times G$, there exists a homomorphism
$\rho: G\times G \to G$
such that $\rho(h_i)=1$ if and only if $h_i=1$.
\end{df}

\begin{lm}
\label{disc}
Let $G$ be a discriminating group then for any integers $n$, $N$
and any elements
$h_1,\dots,h_N \in G^{\times n}$ there exists a homomorphism
$\rho: G^{\times n}\to G$
such that $\rho(h_i)=1$ if and only if $h_i=1$.
\end{lm}
\begin{proof}
We will use induction on $n$. The base case $n=1$ is trivial (take
$\rho = \id$).

Suppose that we have elements $h_i\in G^{\times n+1}$. We can
express each of them as
$h_i = (h_{i,1},h_{i,2})$, where $h_{i,1} \in G^{\times n}$ and
$h_{i,2} \in G$.
By the induction hypothesis, there exists a homomorphism
$\rho_1: G^{\times n} \to G$,
such that $\rho_1(h_{i,1})=1$ iff $h_{i,1} =1$. Now consider the elements
$\bar h_i = (\rho_i(h_{i,1}),h_{i,2}) \in G \times G$. By the
discriminating property of the group $G$, there exists a map
$\rho_2:G\times G \to G$,
such that $\rho_2(\bar h_i)=1$ iff $\bar h_i=1$. The construction of these maps
show that if we define $\rho = \rho_2 \circ (\rho_1 \times \id)$, then
$\rho (h_i) =1$ iff and only if $\bar h_i =1 $. This is the same as
$ \rho_i(h_{i,1})=1 $ and $h_{i,2}=1$, which is equivalent to $h_i=1$.
This finishes the proof of the induction step and
completes the proof of the lemma.
\end{proof}

\begin{thq}
\label{weakid}
Let $G$ be a discriminating group. If a set $S$ is a set of
weak identities in $G$,
then every element from $S$ is an identity in the group $G$.
\end{thq}
\begin{proof}
Assume the contrary: there exists an $s\in S$ such that $s$ is
not an identity in $G$.
Then there exists a homomorphism $\pi$ from $\mr{F}$ to $G$,
such that $\pi(s)\not=1$.
Let us consider the elements $h_j = i_j(\pi(s)) \in G^{\times N}$
for $j=1,\dots, N$.
By lemma~\ref{disc}, there exists a homomorphism
$\tilde \rho: G^{\times N} \to G$,
such that $\tilde \rho(h_j) \not = 1$ for all $j$. Finally, define
a homomorphism $\rho: \mr{F}^{N} \to G$ by
$$
\rho = \tilde \rho \circ ( \pi \times \pi \times \dots \times \pi).
$$
Then by construction we have that $\rho(i_j(s)) \not = 1$,
which contradicts the
fact that $S$ is a set of weak identities in the group $G$.
\end{proof}

\begin{rem}
The converse
(that if every weak identity is an identity then the group $G$
is discriminating)
is not true --- for example, in any abelian group, all weak
identities are identities
but not all abelian groups are discriminating.
\end{rem}

\begin{co}
\label{weak*id}
Let $G$ be a discriminating group. If a $T$-subgroup $\mr{S}$ is a set
of weak* identities in $G$
modulo the trivial $T$-subgroup,
then every element from $\mr{S}$ is an identity in the group $G$.
\end{co}
\begin{proof}
By the definition of weak* identities, there exists an integer
$n$ and  $T$-subgroups $\mr{S}_i$, for
$i=0,\dots, n$ such that
$\mr{S} = \mr{S}_0$, $\mr{S}_n = \{1\}$,
and $\mr{S}_{i-1} \wId{\mr{S}_i}$, for all $i=1,\dots,n$. By induction on $k$,
we can show that $\mr{S}_{n-k}$ are identities in $G$.
The induction step is done using theorem~\ref{weakid}.
\end{proof}

Apply theorem~\ref{weakid} and corollary~\ref{weak*id} to the results in
sections~\ref{linear} and~\ref{solvable} we obtain the following results:

\begin{thq}
\label{disc:solvable}
A finitely generated solvable group $G$
is discriminating if and only if it is isomorphic to $\Z^n$.
\end{thq}
\begin{proof}
By theorem~\ref{main}, the $T$-subgroup $\mr{F}'$ consists of
weak* identities in
$G$. Now we can apply the previous corollary to show that $\mr{F}'$ consists of
identities in the group $G$, i.e., the group $G$ is abelian.
It is easy to see that the only finitely generated abelian groups,
which are discriminating, are the torsion free ones.
\end{proof}

\begin{thq}
\label{disc:linear}
A linear group $G$ is discriminating only if it is abelian.
\end{thq}

\end{section}


\begin{section}{Interesting questions}
\label{sec:open}

Finally we mention several interesting open questions involving the notion
of weak identities.

\begin{que}
\emph{Is it true that $\mr{wId}(G)$ is always a set of weak
identities in the group $G$. }

As mentioned in remark~\ref{widentities} this is true if
$\mr{wId}(G)$ is finitely generated as a verbal subgroup. This
question is equivalent to the following question: Let $S$ be a
subset of $\mr{F}$ such that any one element subset of $S$ is a
set of weak identities in $G$, is it true that $S$ is a set of
weak identities. The answer is positive if the set $S$ is finite.
\end{que}


\begin{que}
\emph{Is it possible to characterize the classes of groups
$\mr{wVar}(\mr{H})$.}

It can be shown that if a group $G$ is in the class $\mr{wVar}(\mr{H})$ and
$\rho: G\to H$ is surjective homomorphism such that for any finite set of
elements $h_i\in H$ there exist elements $g_i\in G$ such that
$\rho(g_i) = h_i$ and $[g_i,g_j] =1$ whenever $[h_i,h_j]=1$; then the group
$H$ is also in the class $\mr{wVar}(\mr{H})$. This class is closed under
taking subgroups, finite Cartesian products and
`restricted' homomorphic images.
\end{que}


\begin{que}
\label{describe}
\emph{Describe all group in $\mr{wVar}(\mr{F}')$.}

This class contains all groups $G$ which has bounded centralizer sequences, i.e.,
all finite groups, all finitely generated meta-abelian groups, all linear groups
and all free groups.
As mentioned in the previous question, this class is closed under taking subgroups,
finite products and restricted homomorphic images.
Are there examples of groups $G \in \mr{wVar}(\mr{F}')$, which does not have a
bounded centralizer sequence?
\end{que}


\begin{que}
\emph{Describe all finitely generated groups $G$ such that
any weak identity in $G$ is an identity.}

It is easy to see that for any abelian group $G$, every weak
identity in $G$ is an
identity. By theorem~\ref{weakid} the same is true for any
discriminating group. This class is also closed under Cartesian products.
Are there examples of non discriminating, non abelian groups,
which does not decompose as a Cartesian product, with this property?
\end{que}


\begin{que}
\emph{Does there exist a finitely generated group $G$ such that
$$
\mr{Id}(G) \subset \mr{wId}(G) \subset \mr{F}',
$$
where all the inclusions are strict?}

Any finitely generated solvable group $G$
such that $[g_1,g_2]$ is not a weak identity in $G$ will give
positive answer to this question. An interesting generalization is
whether there exists a finitely generated group $G$ such that
$$
\mr{Id}(G) \subset \mr{wId}(G) \subset \mr{w^*Id}(G) \subset \mr{F}',
$$
where all the inclusions are strict.
\end{que}
\end{section}

After completing the work on the paper, the author was informed that
Theorems~\ref{disc:linear} and~\ref{disc:solvable} were proven
independently by A. Myashnikov and P. Shumyatsky~\cite{Ma}
using different methods.


\bibliographystyle{plain}
\bibliography{Disc}

\end{document}